# CORRECTION
# CENTRAL LIMIT THEOREMS FOR ADDITIVE FUNCTIONALS OF THE SIMPLE EXCLUSION PROCESS

BY S. SETHURAMAN

*Iowa State University*

Definition 2.1 in the above paper is incorrectly stated. In the proof of Theorem 2.1, which gives an invariance principle for certain processes satisfying Definition 2.1, conditions in Definition 2.1 are sufficient to deduce finite-dimensional convergence, but not enough to apply a maximal inequality for "demimartingales" to obtain tightness. The problem is Definition 2.1, as stated, only considers "pair increment associations" and not more general associations needed for the demimartingale property. We slightly strengthen the definition here in this correction so that the proof of tightness in Theorem 2.1 holds. Details of how this is accomplished are given below.

By substituting the corrected Definition 2.1 for the previous one, all results in the article hold as written. In particular, Proposition 2.1, which is the link between Theorem 2.1 and the main results, and which states certain additive processes satisfy Definition 2.1, holds with the same argument.

CORRECTED DEFINITION 2.1. Let $\{\vec{v}(t) = (v_1(t), \ldots, v_m(t)) : t \geq 0\}$ be an $m$-dimensional $L^2$ process with stationary increments. We say $\vec{v}$ has weakly positive associated increments if

$$E[\phi(\vec{v}(t+s) - \vec{v}(s))\psi(\vec{v}(s_1), \ldots, \vec{v}(s_n))]$$
$$\geq E[\phi(\vec{v}(t))]E[\psi(\vec{v}(s_1), \ldots, \vec{v}(s_n))]$$

for all coordinatewise increasing functions $\phi : \mathbb{R}^m \to \mathbb{R}$ and $\psi : (\mathbb{R}^m)^n \to \mathbb{R}$, and all $s, t \geq 0$, $0 \leq s_1 < \cdots < s_n = s$ and $n \geq 1$.

We remark the earlier Definition 2.1 only stipulated the pair condition

$$E[\phi(\vec{v}(t+s) - \vec{v}(s))\psi(\vec{v}(s))] \geq E[\phi(\vec{v}(t))]E[\psi(\vec{v}(s))].$$

---

We now indicate how the modified definition applies in the proof of tightness in Theorem 2.1. Following standard tightness arguments, one needs to prove for a continuous mean-zero scalar process $v(t)$ with stationary increments satisfying corrected Definition 2.1, with $v(0) = 0$, $\lim_{t\to\infty} t^{-1} E[v(t)^2] = \sigma^2$ and $t^{-1/2} v(t) \Rightarrow N(0, \sigma^2)$ that, for all $\varepsilon > 0$,

$$(1) \qquad \lim_{\delta \downarrow 0} \limsup_{\alpha \to \infty} \frac{1}{\delta} P\left[ \sup_{t \in [0,\delta]} |v(\alpha t)| > \varepsilon \sqrt{\alpha} \right] = 0.$$

For $\delta > 0$, let $A$ be a countable dense set of $[0, \delta]$, and for $n \geq 1$, let $A_n$ be a set of $n$ points so that $A_n \uparrow A$. Fix also that $\delta \in A$ and $\delta \in A_1$. Then, for $\alpha \geq 1$, by continuity $\sup_{t \in [0,\delta]} |v(\alpha t)| = \sup_{t \in A} |v(\alpha t)|$, and for $n$ large enough

$$P\left[\sup_{t \in A} |v(\alpha t)| > \varepsilon \sqrt{\alpha}\right] \leq 2P\left[\sup_{t \in A_n} |v(\alpha t)| > \varepsilon \sqrt{\alpha}\right].$$

Let now $0 \leq t_1 < \cdots < t_{n-1} < t_n = \delta$ be a labeling of $A_n$. From the corrected definition and mean-zero property $E[v(t)] = 0$ we have $E[(v(\alpha t_{j+1}) - v(\alpha t_j))\psi(v(\alpha t_j), \ldots, v(\alpha t_1))] \geq 0$ for all $1 \leq j \leq n-1$ and increasing $\psi$, and so $\{v(\alpha t): t \in A_n\}$ is a demimartingale (cf. [1], page 362). Hence, we can apply the maximal inequality ([2], Corollary 6) and variance convergence $\lim_{\alpha \to \infty} (\alpha \delta)^{-1} E[v(\alpha \delta)^2] = \sigma^2$ to get

$$\limsup_{\alpha \to \infty} P\left[\sup_{t \in A_n} |v(\alpha t)| > \varepsilon \sqrt{\alpha}\right] \leq C_0 \frac{\sigma \sqrt{\delta}}{\varepsilon} \lim_{\alpha \to \infty} \left\{ P\left[|v(\alpha \delta)| > \frac{\varepsilon}{2} \sqrt{\alpha}\right] \right\}^{1/2}$$

for a universal constant $C_0$. From marginal convergence, $\lim_{\alpha \to \infty} P[|v(\alpha \delta)| > (\varepsilon/2)\sqrt{\alpha}] = (2\pi\sigma^2)^{-1/2} \int_{(\varepsilon/2)\delta^{-1/2}}^{\infty} \exp(-x^2/(2\sigma^2)) \, dx$, and so (1) holds.

Also, we note typos: in line 8, page 281, change $2/(\pi \det(\sigma_p^2))$ to $1/(\pi \times (\det(\sigma_p^2))^{1/2})$; in lines 9–10, page 286, $ds$ to $dr$; in line 10, page 293, $=$ to $\geq$; in line 29, page 294, change $> 0$ to $< \infty$; in line 27, page 297, $\exp(\lambda^2 s - \lambda - 1)(-\lambda s)$ should be $(\lambda^2 s - 2\lambda) \exp(-\lambda s)$.

DEPARTMENT OF MATHEMATICS
IOWA STATE UNIVERSITY
396 CARVER HALL
AMES, IOWA 50011
USA
E-MAIL: sethuram@iastate.edu